\documentclass[12pt,oneside]{amsproc}
\usepackage[koi8-r]{inputenc}
\usepackage[russian]{babel}
\usepackage{amssymb,epic,eepic}
\numberwithin{equation}{section}
\newtheorem{lem}{Лемма}[section]
\newtheorem{tm}{Теорема}[section]

\newtheorem{prop}{Утверждение}[section]
\newtheorem{rim}{Замечание}[section]

\title[]{Самоподобные функции в пространстве \(L_2[0,1]\) и задача
Штурма-Лиувилля с сингулярным индефинитным весом}
\author{А.~А.~Владимиров,%
\address{Московский государственный университет
им.~М.~В.~Ломоносова, механико-математический факультет}
\email{vladimi@mech.math.msu.su}
И.~А.~Шейпак%
\address{Московский государственный университет
им.~М.~В.~Ломоносова, механико-математический факультет}
\email{iasheip@mech.math.msu.su}}
\thanks{Работа поддержана РФФИ, грант \No\,04-01-00712, и фондом поддержки
ведущих научных школ, грант~НШ-1927.2003.1.}
\newcommand{\Wo}{{\raisebox{0.2ex}{\(\stackrel{\circ}{W}\)}}{}}
\newcommand{\ind}{\operatorname{ind}}

\begin{document}
\noindent УДК~517.984
\begin{abstract}
В статье изучается вопрос об асимптотике спектра граничной задачи
\begin{gather*}
	-y''-\lambda\rho y=0,\\ y(0)=y(1)=0,
\end{gather*}
где \(\rho\) есть функция из пространства \(\Wo_2^{-1}[0,1]\), имеющая
арифметически самоподобную первообразную. При этом требование
знакоопределённости на вес \(\rho\) не накладывается. Полученные
теоретические результаты иллюстрируются данными численных расчётов.
\end{abstract}
\maketitle

\section{Введение}
Рассмотрим граничную задачу
\begin{gather}\label{eq1:3}
	-y''-\lambda\rho y=0,\\ \label{eq1:4} y(0)=y(1)=0.
\end{gather}
Хорошо известно, что если вес \(\rho\) принадлежит пространству
\(C[0,1]\) и положителен, то граничная задача~\eqref{eq1:3},~\eqref{eq1:4}
эквивалентна спектральной задаче для некоторого неограниченного самосопряжённого
оператора в гильбертовом пространстве \(L_2([0,1];\;\rho)\). Также известно, что
если вес \(\rho\in C[0,1]\) не является знакоопределённым, то
задача~\eqref{eq1:3},~\eqref{eq1:4} эквивалентна спектральной задаче для
неограниченного самосопряжённого оператора в пространстве Крейна
(см.~\cite{CL}).

Ситуация, когда вес \(\rho\) представляет собой сингулярную обобщённую функцию,
изучена существенно хуже. Более того, в настоящее время для такой задачи нет
даже канонической операторной модели. В известных авторам работах по указанной
проблематике~\cite{KacKr} и~\cite{SV} предложены две различные (хотя и
оказывающиеся эквивалентными) интерпретации задачи~\eqref{eq1:3},~\eqref{eq1:4}
с сингулярным весом \(\rho\), притом только для случая, когда \(\rho\)
неотрицателен (то есть представляет собой меру). Случай, когда весовая функция 
не является знакоопределённой, видимо, не исследовался вовсе.

Между тем, разрабатываемая в настоящее время теория операторов Штурма--Лиувилля
с сингулярными потенциалами (см., например,~\cite{SaSh}) позволяет естественным
образом связывать с задачей~\eqref{eq1:3},~\eqref{eq1:4} линейный пучок
ограниченных операторов в том весьма общем случае, когда \(\rho\) есть
произвольный элемент пространства \(\Wo_2^{-1}[0,1]\) (подробности см.~в
параграфе~\ref{pt:3}). Такой пучок может быть построен даже тогда, когда вместо
\(-y''\) в левой части равенства~\eqref{eq1:3} стоит произвольное
полуограниченное дифференциальное выражение второго порядка (хотя в настоящей
работе мы и ограничиваемся рассмотрением выражения \(-y''\)). Для исследования
спектра такого пучка оказывается возможным привлечение вариационной техники из
работы~\cite{LSY} (см.~далее теорему~\ref{tm2:1}). Поэтому принципиальных
трудностей рассмотрение задачи~\eqref{eq1:3},~\eqref{eq1:4} с весовой функцией
\(\rho\in\Wo_2^{-1}[0,1]\) на самом деле не представляет.

Гораздо серьёзнее оказываются трудности технические. В отличие от регулярного
случая \(\rho\in C[0,1]\), в общем случае \(\rho\in\Wo_2^{-1}[0,1]\) даже
порядок асимптотики спектра задачи~\eqref{eq1:3},~\eqref{eq1:4} зависит от
выбора весовой функции (см., например,~\cite{SV}). Поэтому для установления
таких асимптотик необходимо привлекать дополнительную информацию о структуре
весовой функции \(\rho\). Мы будем предполагать, что \(\rho\) представляет
собой обобщённую производную \emph{самоподобной} квадратично суммируемой
функции "--- порядок асимптотики будет тогда определяться параметрами
самоподобия. Таким образом, некоторые результаты настоящей статьи можно
рассматривать как обобщения утверждений из работы~\cite{SV} на более широкий
класс весовых функций.

В работе~\cite{SV} для вывода спектральных асимптотик была использована теория
восстановления. В настоящей статье будет применена эта же техника в её несколько
расширенном варианте (при исследовании некоторых индефинитных задач вместо
обычного уравнения восстановления возникает его "`двумерный"' аналог).

Заслуживает упоминания следующее отличие индефинитного случая от дефинитного.
Для задачи со знакоопределённым весом известна асимптотическая оценка
\begin{equation}\label{eq:KacKr}
	\dfrac{1}{\lambda_n}=O\left(\dfrac{1}{n^2}\right)\qquad\mbox{при }
	n\to\infty
\end{equation}
(см.~\cite[(11.7)]{KacKr}, \cite[(3)]{SV}). Известен также ряд результатов (см.
ссылки в~\cite{SV}), уточняющих эту оценку в зависимости от наличия у обобщённой
первообразной веса \(\rho\) ненулевой абсолютно непрерывной составляющей. Для
индефинитного случая соответствующая проблематика теряет смысл, так как
оценка~\eqref{eq:KacKr} опровергается на примерах. Действительно,
оценка~\eqref{eq:KacKr} имеет место для тех и только тех весовых функций, у
которых спектральный порядок обобщённой первообразной (см.
параграф~\ref{pt:666}) не превосходит \(1\). Однако для произвольной
(немонотонной) самоподобной функции из пространства \(L_2[0,1]\) спектральным
порядком может быть любое число, меньшее чем \(2\) (см. лемму~\ref{lem:D} и
утверждение~\ref{prop4:1}).

Структура статьи такова. В параграфе~\ref{pt:2} даётся доказательство одного
частного случая теоремы восстановления и её двумерного аналога. При этом
используется модификация метода из работы~\cite{LV}. В параграфе~\ref{pt:666}
очерчиваются контуры теории самоподобных функций в пространстве \(L_2[0,1]\). В
параграфе~\ref{pt:3} исследуются спектральные свойства отвечающего
задаче~\eqref{eq1:3},~\eqref{eq1:4} операторного пучка. Для случая, когда
весовая функция \(\rho\) оказывается обобщённой производной арифметически
самоподобной функции положительного спектрального порядка, выписываются
асимптотики спектра (случай неарифметического самоподобия будет разобран в
следующей статье). Наконец, в параграфе~\ref{pt:4} полученные результаты
иллюстрируются на примерах, в том числе приводятся данные численных расчётов.

\section{Теоремы восстановления}\label{pt:2}
Введём в рассмотрение банаховы пространства \(\ell_{1,r}\), где
\(r\in\mathbb R^+\), элементами которых являются последовательности
\(\theta=\{\theta_k\}_{k=0}^{\infty}\), ограниченные по норме
\[
	\|\theta\|_{\ell_{1,r}}=\sum\limits_{k=0}^{\infty}r^k\,|\theta_k|.
\]
Произвольному элементу \(\theta\in\ell_{1,r}\) можно сопоставить определённую
на круге \(\{w\;\vline\;|w|<r\}\) аналитическую функцию \(\Theta\) вида
\[
	\Theta(w)=\sum\limits_{k=0}^{\infty}\theta_k w^k.
\]
В дальнейшем \(\Theta\) будет называться \emph{производящей функцией}
вектора \(\theta\).

\begin{lem}\label{lemA:0}
Пусть фиксированы произвольные \(r<1\) и \(R>1\). Пусть также набор
неотрицательных чисел \(\{u_{k}\}_{k=1}^{N}\) таков, что
\[
	\sum\limits_{k=1}^N u_{k}=1,
\]
причём наибольший общий делитель номеров \(k\), для которых справедливо
неравенство \(u_k>0\), равен \(1\). Тогда для любого \(x\in\ell_{1,R}\)
существует и единственно решение \(z\in\ell_{1,r}\) системы
\[
	\forall n\in\{0,1,2,\ldots\}\qquad z_{n}=x_{n}+
	\sum\limits_{k=1}^{\min(N,n)}u_kz_{n-k}.
\]
При этом координаты вектора \(z\) удовлетворяют оценке
\[
	\forall n\in\{0,1,2,\ldots\}\qquad \left|z_{n}-
	\dfrac{\omega}{J}\right|\leqslant\|x\|_{\ell_{1,R}}\cdot C_{n},
\]
где бесконечно малая последовательность \(\{C_k\}_{k=0}^{\infty}\) не
зависит от выбора \(x\), а величины \(\omega\) и \(J\) определены равенствами
\begin{gather*}
	\omega=\sum\limits_{k=0}^{\infty} x_{k}\\
	\intertext{и}
	J=\sum\limits_{k=1}^N ku_k.
\end{gather*}
\end{lem}

Для доказательства леммы~\ref{lemA:0} достаточно воспроизвести, с
незначительными изменениями, рассуждения из доказательства
леммы~\cite[Lemma~A.1]{LV}.

\begin{lem}\label{lemA:1}
Пусть фиксированы произвольные \(r<1\) и \(R>1\). Пусть наборы неотрицательных
чисел \(\{u_{k}\}_{k=1}^{N}\) и \(\{v_{k}\}_{k=1}^{N}\) таковы, что
\[
	\sum\limits_{k=1}^N (u_{k}+v_{k})=1,\qquad
	\sum\limits_{k=1}^N v_{k}>0,
\]
причём наибольший общий делитель номеров \(k\), для которых справедливо
неравенство \(u_k+v_k>0\), равен \(1\). Пусть также найдётся либо нечётное
\(k\leqslant N\) со свойством \(u_k>0\), либо чётное \(k\leqslant N\) со
свойством \(v_k>0\). Тогда для любых \(x_j\in\ell_{1,R}\), где \(j=1,2\),
существует и единственна пара \(z_j\in\ell_{1,r}\) решений системы
\[
	\forall n\in\{0,1,2,\ldots\}\qquad z_{j,n}=x_{j,n}+
	\sum\limits_{k=1}^{\min(N,n)}(u_kz_{j,n-k}+v_kz_{3-j,n-k}),
	\qquad j=1,2.
\]
При этом координаты векторов \(z_j\), где \(j=1,2\), удовлетворяют оценке
\begin{equation}\label{eq4:3}
	\forall n\in\{0,1,2,\ldots\}\qquad\left|z_{j,n}-\dfrac{\omega}{J}\right|
	\leqslant(\|x_1\|_{\ell_{1,R}}+\|x_2\|_{\ell_{1,R}})\cdot C_{n},
	\qquad j=1,2,
\end{equation}
где бесконечно малая последовательность \(\{C_k\}_{k=0}^{\infty}\) не
зависит от выбора \(x_j\), а величины \(\omega\) и \(J\) определены равенствами
\begin{gather}\label{eqA:-1}
	\omega=\dfrac12\sum\limits_{k=0}^{\infty}(x_{1,k}+x_{2,k})\\
	\intertext{и} \label{eqA:1}
	J=\sum\limits_{k=1}^N k(u_k+v_k).
\end{gather}
\end{lem}
\begin{proof}
\textit{Шаг~1.}
Обозначим через \(U\) и \(V\) многочлены
\begin{align*}
	U(w)&=\sum\limits_{k=1}^Nu_kw^k,\\
	V(w)&=\sum\limits_{k=1}^Nv_kw^k.
\end{align*}
Из условий леммы вытекает, что многочлен \(1-U-V\) имеет в круге
\(\{w\mid |w|\leqslant 1\}\) единственный, причём простой, нуль
\(w=1\). Тем самым рассматриваемый многочлен может быть представлен
в виде \((1-w)\cdot Q(w)\), где \(Q\) есть многочлен вида
\[
	Q(w)=\sum\limits_{n=0}^{N-1}\left(\sum\limits_{k=n+1}^N(u_k+v_k)\right)w^n.
\]
В свою очередь, \(Q\) может быть представлен в виде \((1-w)\cdot P(w)+J\), где
\(P\) есть некоторый многочлен, а постоянная \(J=Q(1)\) определена
равенством~\eqref{eqA:1}.

Из условий леммы вытекает также, что при любом \(w\), удовлетворяющем
неравенству \(|w|\leqslant 1\), справедливо неравенство
\[
	|U(w)-V(w)|<1.
\]
Объединяя этот факт со сказанным ранее, получаем, что при некотором
\(\hat R\in(1,R)\) круг \(\{w\mid |w|<\hat R\}\) свободен от нулей многочленов
\(Q\) и \(1-U+V\).

\textit{Шаг~2.}
Производящие функции \(X_j\) и \(Z_j\), где \(j=1,2\), векторов \(x_j\) и
\(z_j\) должны подчиняться равенствам
\[
	Z_j=X_j+UZ_j+VZ_{3-j},\qquad j=1,2,
\]
из которых немедленно вытекают равенства
\begin{equation}\label{eqA:3}
	Z_j=\dfrac{X_1+X_2}{2(1-U-V)}+\dfrac{X_j-X_{3-j}}{2(1-U+V)},\qquad j=1,2.
\end{equation}
Последние равенства доказывают факт существования и единственности пары решений
\(z_j\), где \(j=1,2\).

Далее, равенства~\eqref{eqA:3} означают, что при любом \(w\) из области
определения функций \(Z_j\) справедливы равенства
\begin{multline*}
	Z_j(w)=\dfrac{X_1(1)+X_2(1)}{2J(1-w)}+\left[\dfrac{(X_1(w)-X_1(1))+
	(X_2(w)-X_2(1))}{2J(1-w)}-\right.\\ \left.-\dfrac{P(w)\cdot
	(X_1(w)+X_2(w))}{2JQ(w)}+\dfrac{X_j(w)-X_{3-j}(w)}{2(1-U(w)+V(w))}
	\right].
\end{multline*}
Слагаемое из правой части, заключённое в квадратные скобки, представляет
собой аналитическую функцию в круге \(\{w\;\vline\;|w|<\hat R\}\).
Применяя к этой функции неравенства Коши, убеждаемся, что её
коэффициенты Маклорена \(c_k\), где \(k=0,1,2,\ldots\), удовлетворяют
оценке
\[
	\forall k\in\{0,1,2,\ldots\}\qquad |c_k|\leqslant C_k
	(\|x_1\|_{\ell_{1,R}}+\|x_2\|_{\ell_{1,R}}),
\]
где бесконечно малая последовательность \(\{C_k\}_{k=0}^{\infty}\)
не зависит от выбора \(x_j\). Поскольку коэффициенты Маклорена функции
\[
	\dfrac{X_1(1)+X_2(1)}{2J(1-w)}
\]
тождественно равны \(\omega/J\), где постоянная \(\omega\) определена
равенством~\eqref{eqA:-1}, то тем самым неравенства~\eqref{eq4:3}
справедливы. Лемма доказана.
\end{proof}

Леммы~\ref{lemA:0} и~\ref{lemA:1} позволяют доказать следующие теоремы
восстановления с оценкой остатка.

\begin{tm}\label{tmA:1}
Пусть фиксировано произвольное положительное вещественное число \(\tau\).
Пусть также набор неотрицательных чисел \(\{u_{k}\}_{k=1}^{N}\) таков, что
\[
	\sum\limits_{k=1}^N u_{k}=1,
\]
причём наибольший общий делитель номеров \(k\), для которых справедливо
неравенство \(u_k>0\), равен \(1\). Наконец, пусть непрерывная на
\(\mathbb R\) функция \(X\) удовлетворяет условиям
\begin{align*}
	&\forall t\in\mathbb R^+\qquad |X(t)|\leqslant \Pi\,e^{-\tau t},\\
	&\forall t\in\mathbb R^-\qquad X(t)=0,
\end{align*}
где \(\Pi>0\). Тогда существует и единственна непрерывная на
\(\mathbb R\) функция \(Z\), удовлетворяющая уравнениям
\begin{align*}
	&\forall t\in\mathbb R^+\qquad Z(t)=X(t)+\sum\limits_{k=1}^N u_k
	\,Z(t-k),\\
	&\forall t\in\mathbb R^-\qquad Z(t)=0.
\end{align*}
При этом для функции \(Z\) справедлива оценка
\[
	\forall t\in\mathbb R^+\qquad \left|Z(t)-s(t)\right|\leqslant
	\Pi\cdot C(t),
\]
где исчезающая на бесконечности функция \(C\) не зависит от выбора
функции \(X\), а непрерывная \mbox{\(1\)-пе}\-ри\-о\-ди\-чес\-кая
функция \(s\) имеет вид
\[
	\forall t\in\mathbb R\qquad s(t)=\dfrac{1}{J}
	\sum\limits_{k=-\infty}^{+\infty} X(t-k).
\]
Через \(J\) здесь обозначена величина
\[
	J=\sum\limits_{k=1}^N k\,u_k.
\]
\end{tm}
\begin{tm}\label{tmA:2}
Пусть фиксировано произвольное положительное вещественное число \(\tau\).
Пусть наборы неотрицательных чисел \(\{u_{k}\}_{k=1}^{N}\) и
\(\{v_{k}\}_{k=1}^{N}\) таковы, что
\[
	\sum\limits_{k=1}^N (u_{k}+v_{k})=1,\qquad
	\sum\limits_{k=1}^N v_{k}>0,
\]
причём наибольший общий делитель номеров \(k\), для которых справедливо
неравенство \(u_k+v_k>0\), равен \(1\). Пусть также найдётся либо нечётное
\(k\leqslant N\) со свойством \(u_k>0\), либо чётное \(k\leqslant N\) со
свойством \(v_k>0\). Наконец, пусть непрерывные на \(\mathbb R\) функции
\(X_j\), где \(j=1,2\), удовлетворяют условиям
\begin{align*}
	&\forall t\in\mathbb R^+\qquad |X_j(t)|\leqslant \Pi_j\,e^{-\tau t},\\
	&\forall t\in\mathbb R^-\qquad X_j(t)=0.
\end{align*}
Здесь \(\Pi_j\), где \(j=1,2\), "--- положительные числа. Тогда существует и
единственна пара непрерывных на \(\mathbb R\) функций \(Z_j\), где \(j=1,2\),
удовлетворяющая уравнениям
\begin{align*}
	&\forall t\in\mathbb R^+\qquad Z_j(t)=X_j(t)+\sum\limits_{k=1}^N (u_k
	\,Z_j(t-k)+v_k\,Z_{3-j}(t-k)),\\
	&\forall t\in\mathbb R^-\qquad Z_j(t)=0.
\end{align*}
При этом для функций \(Z_j\), где \(j=1,2\), справедливы оценки
\[
	\forall t\in\mathbb R^+\qquad \left|Z_j(t)-s(t)\right|\leqslant
	(\Pi_1+\Pi_2)\cdot C(t),
\]
где исчезающая на бесконечности функция \(C\) не зависит от выбора
функций \(X_j\), где \(j=1,2\), а непрерывная \mbox{\(1\)-пе}\-ри\-о\-ди\-чес\-кая
функция \(s\) имеет вид
\[
	\forall t\in\mathbb R\qquad s(t)=\dfrac{1}{2J}
	\sum\limits_{k=-\infty}^{+\infty}\left(X_1(t-k)+X_2(t-k)\right).
\]
Через \(J\) здесь обозначена величина
\[
	J=\sum\limits_{k=1}^N k\,(u_k+v_k).
\]
\end{tm}

\begin{proof}[Доказательство теоремы~\ref{tmA:2}]
Зафиксируем произвольное \(R\in(1,e^{\tau})\). Тогда при любом
\(t\in[0,1)\) последовательности \(x_j(t)=\{x_{j,k}(t)\}_{k=0}^{\infty}\),
где \(j=1,2\), вида
\[
	\forall k\in\{0,1,2,\ldots\}\qquad x_{j,k}(t)=X_j(t+k),\qquad j=1,2,
\]
принадлежат пространству \(\ell_{1,R}\). При этом равномерно по \(t\in[0,1)\)
справедливы оценки
\[
	\|x_j(t)\|_{\ell_{1,R}}\leqslant\dfrac{\Pi_j}{1-Re^{-\tau}},\qquad j=1,2.
\]

Пусть теперь \(\{C_k\}_{k=0}^{\infty}\) "--- последовательность из
леммы~\ref{lemA:1}. Зафиксируем исчезающую на бесконечности непрерывную
функцию \(C\), удовлетворяющую условию
\[
	\forall t\in\mathbb R^+\qquad C(t)>\dfrac{C_{[t]}}{1-Re^{-\tau}},
\]
где через \([t]\) обозначена целая часть числа \(t\). Утверждение
теоремы теперь легко получается в результате применения леммы~\ref{lemA:1}
к последовательностям \(z_j(t)=\{z_{j,k}(t)\}_{k=0}^{\infty}\), где \(j=1,2\), вида
\[
	\forall k\in\{0,1,2,\ldots\}\qquad z_{j,k}(t)=Z_j(t+k),
\]
где \(t\in[0,1)\).
\end{proof}

Теорема~\ref{tmA:1} доказывается аналогично теореме~\ref{tmA:2} с тем отличием,
что вместо леммы~\ref{lemA:1} используется лемма~\ref{lemA:0}.

\section{Самоподобные функции в пространстве \(L_2[0,1]\)}\label{pt:666}
\subsection{Операторы подобия в пространстве \(L_2[0,1]\)}
Пусть фиксировано натуральное число \(n>1\), и пусть вещественные числа
\(a_k>0\), \(d_k\) и \(\beta_k\), где \(k=1,\ldots,n\), таковы, что
\[
	\sum\limits_{k=1}^n a_k=1.
\]
Данному набору чисел можно поставить в соответствие непрерывный
нелинейный оператор \(G:L_2[0,1]\to L_2[0,1]\) вида
\begin{equation}\label{eq:auxto}
	G(f)=\sum\limits_{k=1}^n\left\{\beta_k\cdot\chi_{(\alpha_k,\alpha_{k+1})}
	+d_k\cdot G_k(f)\right\},
\end{equation}
где использованы следующие обозначения:
\begin{enumerate}
\item через \(\alpha_k\), где \(k=1,2,\ldots,n+1\), обозначены числа \(\alpha_1=0\)
и \(\alpha_k=\sum_{l=1}^{k-1}a_l\), где \(k=2,\ldots,n+1\);
\item через \(\chi_{\Gamma}\), где \(\Gamma\) "--- интервал, обозначена характеристическая
функция интервала \(\Gamma\), рассматриваемая как элемент пространства \(L_2[0,1]\);
\item через \(G_k\), где \(k=1,\ldots,n\), обозначены непрерывные линейные операторы
в пространстве \(L_2[0,1]\), действующие на характеристическую функцию
\(\chi_{(\zeta,\xi)}\) произвольного интервала \((\zeta,\xi)\subset[0,1]\) согласно
правилу
\begin{equation}\label{eq3:2.1}
	G_k(\chi_{(\zeta,\xi)})=\chi_{(\alpha_k+a_k\zeta,\alpha_k+a_k\xi)},
	\qquad k=1,\ldots,n.
\end{equation}
Это правило, как нетрудно видеть, определяет операторы \(G_k\) однозначно.
\end{enumerate}

Операторы \(G\) вида~\eqref{eq:auxto} будут называться \emph{операторами подобия}.

\begin{lem}\label{lem3:1}
Оператор подобия \(G\) является сжимающим в том и только том случае, когда справедливо
неравенство
\begin{equation}\label{eq3:-+}
	\sum\limits_{k=1}^n a_k\,|d_k|^2<1.
\end{equation}
\end{lem}
\begin{proof}
Утверждение доказываемой леммы следует из того факта, что при любых
\(f_1,f_2\in L_2[0,1]\) справедливы  соотношения
\begin{multline*}
	\|G(f_1)-G(f_2)\|_{L_2[0,1]}^2=\int\limits_0^1|G(f_1)-G(f_2)|^2\,dx=\\=
	\sum\limits_{k=1}^n\left(|d_k|^2\,\int\limits_{\alpha_k}^{\alpha_{k+1}}
	|G_k(f_1)-G_k(f_2)|^2\,dx\right)=\left(\sum\limits_{k=1}^n a_k\,|d_k|^2
	\right)\,\int\limits_0^1|f_1-f_2|^2\,dx=\\=\left(\sum\limits_{k=1}^n a_k\,
	|d_k|^2\right)\,\|f_1-f_2\|_{L_2[0,1]}^2.
\end{multline*}
\end{proof}

Из леммы~\ref{lem3:1} и принципа сжимающих отображений немедленно
вытекает справедливость следующего утверждения.
\begin{tm}\label{sek3:1}
Если справедливо неравенство~\eqref{eq3:-+}, то существует и единственна
функция \(f\in L_2[0,1]\), удовлетворяющая уравнению \(G(f)=f\).
\end{tm}

В дальнейшем всегда будет предполагаться, что неравенство~\eqref{eq3:-+}
выполнено.

\subsection{Самоподобные функции и их типы}
Если функция \(f\in L_2[0,1]\) удовлетворяет уравнению \(G(f)=f\), где
\(G\) "--- некоторый оператор подобия, то такая функция будет называться
\emph{самоподобной}. При этом величины \(n\), \(a_k\), \(d_k\) и \(\beta_k\), где
\(k=1,2,\ldots,n\), определяющие соответствующий оператор подобия \(G\), будут
называться \emph{параметрами самоподобия} функции \(f\).

\begin{rim}\label{rim3:1}
Различные операторы подобия могут определять одну и ту же самоподобную функцию. В качестве
примера можно привести функцию \(f\) вида \(f(x)=x\), обладающую, в частности, параметрами
самоподобия
\begin{enumerate}
\item \(n=2\), \(a_1=a_2=d_1=d_2=1/2\), \(\beta_1=0\), \(\beta_2=1/2\);
\item \(n=2\), \(a_1=d_1=(3-\sqrt{5})/2\), \(a_2=d_2=(\sqrt{5}-1)/2\),
\(\beta_1=0\), \(\beta_2=(3-\sqrt{5})/2\);
\item \(n=2\), \(a_1=d_1=1/3\), \(a_2=d_2=2/3\), \(\beta_1=0\), \(\beta_2=1/3\).
\end{enumerate}
\end{rim}

Функция \(f\in L_2[0,1]\), для которой найдутся такие параметры самоподобия, что при
некотором \(\nu>0\) будет справедливо условие
\begin{gather}\label{eq:l_k666}
	\forall k\in\{1,\ldots,n\}\;\exists l_k\in\mathbb N\qquad
	(a_k\,|d_k|)\cdot(a_k\,|d_k|-e^{-l_k\,\nu})=0,
\end{gather}
будет называться \emph{арифметически самоподобной} функцией. Если для некоторых параметров
самоподобия число \(\hat\nu\) является максимальным среди чисел \(\nu\) со
свойством~\eqref{eq:l_k666}, то такое число \(\hat\nu\) будет называться \emph{шагом}
самоподобия функции \(f\).

Функция \(f\in L_2[0,1]\), для которой найдутся такие параметры самоподобия, что при
любом \(\nu>0\) условие~\eqref{eq:l_k666} будет нарушено, будет называться
\emph{неарифметически самоподобной} функцией.

\begin{rim}
Одна и та же функция может иметь различные шаги самоподобия, и даже может быть арифметически
и неарифметически самоподобной одновременно. Например, упоминавшиеся в замечании~\ref{rim3:1}
наборы параметров самоподобия для функции \(f\) вида \(f(x)=x\) показывают, что эта функция:
\begin{enumerate}
\item арифметически самоподобна с шагом \(\ln 4\);
\item арифметически самоподобна с шагом \(\ln 2-\ln(3-\sqrt{5})\);
\item неарифметически самоподобна.
\end{enumerate}
\end{rim}

\subsection{Спектральный порядок} Особое место среди самоподобных функций занимают
функции, для которых существуют параметры самоподобия со следующими свойствами:
\begin{enumerate}
\item среди чисел \(d_k\), где \(k=1,\ldots,n\), не менее двух отличны от нуля;
\item среди чисел \(\beta_k\), где \(k=1,\ldots,n\), по меньшей мере одно отлично
от нуля.
\end{enumerate}
Такие самоподобные функции будут называться \emph{самоподобными функциями положительного
спектрального порядка}.

\begin{lem}\label{lem:D}
Пусть \(f\) "--- самоподобная функция, и пусть \(n\), \(a_k\) и \(d_k\), где
\(k=1,\ldots,n\), "--- её параметры самоподобия. Пусть при этом среди чисел \(d_k\) не менее
двух отличны от нуля. Тогда существует и единственно положительное решение \(D\) уравнения
\begin{equation}\label{eq:spord}
	\sum\limits_{k=1}^n\left(a_k\,|d_k|\right)^{D/2}=1.
\end{equation}
При этом \(D<2\).
\end{lem}
\begin{proof}
Рассмотрим определённую на \((0,+\infty)\) функцию \(\Upsilon\) вида
\[
	\forall x\in(0,+\infty)\qquad\Upsilon(x)=\sum\limits_{k=1}^n\left(a_k\,|d_k|\right)^x.
\]
Из легко получаемых на основе неравенства Коши--Буняковского соотношений
\begin{equation}\label{eq666:666}
	\sum\limits_{k=1}^n a_k\,|d_k|\leqslant\left(\sum\limits_{k=1}^n a_k\,|d_k|^2\right)^{1/2}
	\cdot\left(\sum\limits_{k=1}^n a_k\cdot 1^2\right)^{1/2}=\left(\sum\limits_{k=1}^n a_k\,
	|d_k|^2\right)^{1/2}<1
\end{equation}
следует, что эта функция является убывающей, причём для любого \(x\geqslant 1\) справедливо
неравенство \(\Upsilon(x)<1\). С другой стороны, из условий леммы следует, что при
\(x\approx 0\) значения \(\Upsilon(x)\) превосходят \(1\). Тем самым все утверждения
леммы справедливы.
\end{proof}

При помощи рассуждений, аналогичных тем, что будут ниже проделаны при доказательстве
теоремы~\ref{tm3:2}, можно установить справедливость следующего утверждения.
\begin{tm}\label{tm:spord}
Пусть \(f\in L_2[0,1]\) "--- самоподобная функция положительного спектрального порядка.
Тогда решение \(D\) уравнения~\eqref{eq:spord} не зависит от выбора параметров
самоподобия.
\end{tm}

Решение \(D\) уравнения~\eqref{eq:spord}, отвечающего самоподобной функции положительного
спектрального порядка, будет называться \emph{спектральным порядком} этой функции.
Спектральным порядком самоподобной функции, не являющейся самоподобной функцией
положительного спектрального порядка, будет по определению считаться число \(0\).

\section{Операторный пучок и его спектральные свойства}\label{pt:3}
\subsection{Построение операторной модели}
В дальнейшем через \(\mathcal H\) будет обозначаться пространство
\(\Wo_2^1[0,1]\), снабжённое нормой
\[
	\forall y\in\Wo_2^1[0,1]\qquad \|y\|_{\mathcal H}=\|y'\|_{L_2[0,1]}.
\]
Через \(\mathcal H'=\Wo_2^{-1}[0,1]\) в дальнейшем будет обозначаться
пространство, дуальное к \(\mathcal H\) относительно \(L_2[0,1]\), то есть
пополнение пространства \(L_2[0,1]\) по норме
\[
	\forall y\in L_2[0,1]\qquad \|y\|_{\mathcal H'}=
	\sup\limits_{\|z\|_{\mathcal H}=1}|\langle y,z\rangle_{L_2[0,1]}|.
\]
Через \(\mathfrak I[y,z]\), где \(y\in\mathcal H'\) и \(z\in\mathcal H\),
будет обозначаться полуторалинейная форма, являющаяся продолжением по непрерывности
формы
\[
	\forall y\in L_2[0,1]\;\forall z\in\mathcal H\qquad
	\mathfrak I[y,z]=\int\limits_0^1 y\overline{z}\,dx.
\]

Как несложно проверить, любой функции \(P\in L_2[0,1]\) можно поставить в соответствие
однозначно определённую функцию \(\rho\in\mathcal H'\) со свойством
\[
	\forall y\in\mathcal H\qquad \mathfrak I[\rho,y]=
	-\int\limits_0^1 P\overline{y'}\,dx.
\]
Такая функция \(\rho\) будет называться \emph{производной} от функции \(P\). Легко
проследить связь введённого определения с известным в теории обобщённых функций
понятием обобщённой производной.

В соответствии с мультипликаторной трактовкой задач Штурма--Лиувилля с
сингулярными потенциалами (см., например, \cite{SaSh}), выберем в качестве
операторной модели для задачи~\eqref{eq1:3},~\eqref{eq1:4} линейный пучок
\(T_{\rho}:\mathcal H\to\mathcal H'\) ограниченных операторов, удовлетворяющий тождеству
\begin{gather}\notag
	\forall\lambda\in\mathbb C\;\forall y,z\in\mathcal H\qquad
	\mathfrak I[T_{\rho}(\lambda)y,z]=\int\limits_0^1 y'\overline{z'}\,dx
	-\lambda\cdot\mathfrak I[\rho,\overline{y}z],\\
	\intertext{В случае, когда вес \(\rho\) представляет собой производную функции
	\(P\in L_2[0,1]\), последнее тождество переписывается в виде}\label{eq2:-3}
	\forall\lambda\in\mathbb C\;\forall y,z\in\mathcal H\qquad
	\mathfrak I[T_{\rho}(\lambda)y,z]=\int\limits_0^1 \bigl\{y'\overline{z'}+
	\lambda P\cdot(y'\overline{z}+y\overline{z'})\bigr\}\,dx.
\end{gather}
Несложно убедиться, что в регулярном случае \(\rho\in C[0,1]\) уравнение
\(T_{\rho}(\lambda)y=0\) эквивалентно задаче~\eqref{eq1:3},~\eqref{eq1:4},
понимаемой обычным образом. Очевидна также справедливость тождества
\begin{equation}\label{eq2:444}
	\forall y\in\mathcal H\qquad \mathfrak I[T_{\rho}(0)y,y]
	=\|y\|^2_{\mathcal H}.
\end{equation}

Прежде чем сформулировать основной результат настоящего пункта, введём
следующие понятия.
\begin{itemize}
\item \emph{Индексом инерции} \(\ind S\) оператора \(S:\mathcal H\to
\mathcal H'\) будет называться максимум размерностей подпространств
\(\mathcal M\subset\mathcal H\), удовлетворяющих условию
\[
	\exists\varepsilon>0\;\forall y\in\mathcal M\quad
	\mathfrak I[Sy,y]\leqslant-\varepsilon\|y\|^2_{\mathcal H}.
\]
\item Собственное значение \(\lambda\in\mathbb C\) линейного пучка
\(S:\mathcal H\to\mathcal H'\) вида \(S(\lambda)=A-\lambda B\) будет
называться \emph{собственным значением положительного (отрицательного) типа},
если справедливо неравенство
\[
	\exists\varepsilon>0\;\forall y\in\ker S(\lambda)\quad
	\mathfrak I[By,y]\geqslant\varepsilon\|y\|^2_{\mathcal H}
\]
(соответственно,
\[
	\left.\exists\varepsilon>0\;\forall y\in\ker S(\lambda)\quad
	\mathfrak I[By,y]\leqslant-\varepsilon\|y\|^2_{\mathcal H}\right).
\]
\end{itemize}

Основной результат настоящего пункта состоит в следующем.
\begin{tm}\label{tm2:1}
Спектр пучка \(T_{\rho}\) чисто дискретен, и все его собственные значения являются простыми.

Все собственные значения пучка \(T_{\rho}\), расположенные правее нуля,
имеют положительный тип, а все собственные значения пучка \(T_{\rho}\),
расположенные левее нуля, имеют отрицательный тип. Для любого \(\lambda>0\) число
собственных значений пучка \(T_{\rho}\), принадлежащих интервалу
\((0,\lambda)\), совпадает с индексом инерции \(\ind T_{\rho}(\lambda)\)
оператора \(T_{\rho}(\lambda)\). Аналогично, для любого \(\lambda<0\) число
собственных значений пучка \(T_{\rho}\), принадлежащих интервалу \((\lambda,0)\),
совпадает с индексом инерции \(\ind T_{\rho}(\lambda)\) оператора \(T_{\rho}(\lambda)\).
\end{tm}
\begin{proof}
Производная \(T_{\rho}'\) пучка \(T_{\rho}\) представляет собой компактный оператор.
Поэтому к пучку \(T_{\rho}\) могут быть приложены результаты работы~\cite{LSY}.

Тождество~\eqref{eq2:444} означает, что пучок \(T_{\rho}\) является сильно
дефинизируемым (см.~теорему~\cite[Theorem~1]{LSY}). Утверждения теоремы
вытекают потому из утверждения~\cite[Proposition~6]{LSY},
теоремы~\cite[Theorem~1]{LSY}, тождества~\eqref{eq2:444} и простого свойства
\[
	\forall\lambda\in\mathbb C\qquad\dim\ker T_{\rho}(\lambda)\leqslant 1.
\]
\end{proof}

Теорема~\ref{tm2:1} показывает, что знание асимптотики при \(\lambda\to\pm\infty\)
для индекса инерции операторов \(T_{\rho}(\lambda)\) эквивалентно знанию
асимптотики при \(\lambda\to\pm\infty\) для собственных значений пучка \(T_{\rho}\).
Поэтому получаемые ниже утверждения о поведении величины \(\ind T_{\rho}(\lambda)\)
легко могут быть переформулированы в утверждения о спектре пучка \(T_{\rho}\).

\subsection{Асимптотики спектра}
Основным результатом настоящего пункта является следующее утверждение.

\begin{tm}\label{tm3:2}
Пусть \(P\in L_2[0,1]\) "--- арифметически самоподобная функция с шагом
\(\nu\), имеющая положительный спектральный порядок \(D\). Пусть при этом
найдётся номер \(k\leqslant n\), для которого выполнено одно из следующих
условий:
\begin{itemize}
\item Справедливо неравенство \(d_k>0\), и отношение
\begin{equation}\label{eq:rel}
	\dfrac{\ln(a_k\,|d_k|)}{\nu}
\end{equation}
нечётно.
\item Справедливо неравенство \(d_k<0\), и отношение~\eqref{eq:rel} чётно.
\end{itemize}
Тогда для пучка~\eqref{eq2:-3} справедливы следующие утверждения.
\begin{enumerate}
\item Существуют такие непрерывные неотрицательные \mbox{\(1\)-пе}\-ри\-оди\-чес\-кие функции
\(s_{\pm}\), что при \(\lambda\to\pm\infty\) справедливы асимптотические представления
\begin{equation}\label{eq3:as}
	\ind T_{\rho}(\lambda)=|\lambda|^{D/2}\,\left(s_{\pm}\left(
	\dfrac{\ln|\lambda|}{\nu}\right)+o(1)\right).
\end{equation}
\item\label{ok:2} Если при некотором \(k\in\{1,\ldots,n\}\) имеет место неравенство \(d_k<0\),
то справедливо тождество
\begin{equation}\label{eq3:id+}
	\forall t\in\mathbb R\qquad s_+(t)=s_-(t).
\end{equation}
\item\label{ok:3} Если для некоторой функции \(y\in\mathcal H\) имеет место неравенство
\(\mathfrak I[\rho,|y|^2]>0\), то функция \(s_+\) положительна. Аналогично,
если для некоторой функции \(y\in\mathcal H\) имеет место неравенство
\(\mathfrak I[\rho,|y|^2]<0\), то функция \(s_-\) положительна.
\end{enumerate}
\end{tm}

Доказательство теоремы~\ref{tm3:2} будет опираться на следующую лемму.

\begin{lem}\label{lem3:2.1}
Пусть функция \(P\in L_2[0,1]\) самоподобна, и пусть \(n\), \(a_k\) и \(d_k\),
где \(k=1,2,\ldots,n\), "--- её параметры самоподобия. Тогда для пучка~\eqref{eq2:-3}
справедливы неравенства
\begin{equation}\label{eq3:3.2}
	\forall\lambda\in\mathbb R\qquad 0\leqslant\ind T_{\rho}(\lambda)-
	\sum\limits_{k=1}^n\ind T_{\rho}(a_k d_k\cdot\lambda)\leqslant n-1.
\end{equation}
\end{lem}
\begin{proof}
Зафиксируем произвольное \(\lambda\in\mathbb R\).

\textit{Шаг~1.} Пусть \(\mathcal M_k\), где \(k=1,2,\ldots,n\) "--- произвольные
подпространства в \(\mathcal H\), удовлетворяющие при некотором \(\varepsilon>0\)
условиям
\[
	\forall y\in\mathcal M_k\qquad\mathfrak I[T_{\rho}(a_k d_k\cdot\lambda)y,y]
	\leqslant-\varepsilon\|y\|^2_{\mathcal H},\qquad k=1,2,\ldots,n.
\]
Рассмотрим в \(\mathcal H\) подпространство
\[
	\mathcal M=\bigoplus\limits_{k=1}^n G_k(\mathcal M_k),
\]
где \(G_k\) "--- операторы, определённые формулой~\eqref{eq3:2.1}. Размерность
подпространства \(\mathcal M\), очевидно, представляет собой сумму размерностей
подпространств \(\mathcal M_k\). С другой стороны, для любого набора функций
\(y_k\in\mathcal M_k\), где \(k=1,2,\ldots,n\), и построенной по нему функции
\(y=\sum_{k=1}^nG_k(y_k)\in\mathcal M\) справедливы соотношения
\begin{multline*}
	\mathfrak I[T_{\rho}(\lambda)y,y]=\sum\limits_{k=1}^n
	\int\limits_{\alpha_k}^{\alpha_{k+1}}\left\{|y'|^2+\lambda P\cdot(y'\overline{y}+
	y\overline{y'})\right\}\,dx=\\=\sum\limits_{k=1}^n\int\limits_0^1
	\left\{a^{-1}_k|y_k'|^2+d_k\cdot\lambda P\cdot(y_k'\overline{y_k}+y_k
	\overline{y_k'})\right\}\,dx=\\=\sum\limits_{k=1}^na^{-1}_k
	\mathfrak I[T_{\rho}(a_k d_k\cdot\lambda)y_k,y_k]\leqslant-\varepsilon
	\sum\limits_{k=1}^n a^{-1}_k\|y_k\|^2_{\mathcal H}=-\varepsilon
	\|y\|^2_{\mathcal H}.
\end{multline*}
Тем самым подпространство \(\mathcal M\) удовлетворяет условию
\begin{equation}\label{eq3:3.x}
	\forall y\in\mathcal M\qquad\mathfrak I[T_{\rho}(\lambda)y,y]\leqslant
	-\varepsilon\|y\|^2_{\mathcal H}.
\end{equation}
Ввиду произвольности выбора подпространств \(\mathcal M_k\), сказанное означает
справедливость левого неравенства в соотношении~\eqref{eq3:3.2}

\textit{Шаг~2.} Пусть \(\mathcal M\) "--- подпространство в \(\mathcal H\), удовлетворяющее
при некотором \(\varepsilon>0\) условию~\eqref{eq3:3.x}. Пусть также \(\mathcal M_k\), где
\(k=1,2,\ldots,n\), "--- подпространства в \(\mathcal H\), удовлетворяющие при некоторых
\(\varepsilon_k>0\) условиям
\begin{equation}\label{eq3:3.xxl}
	\forall y\in\mathcal M_k\qquad\mathfrak I[T_{\rho}(a_k d_k\cdot\lambda)y,y]
	\leqslant-\varepsilon_k\|y\|^2_{\mathcal H},\qquad k=1,2,\ldots,n.
\end{equation}
Через \(\mathcal M_k^{\perp}\) в дальнейшем будут обозначаться подпространства в
\(\mathcal H\) вида
\[
	\mathcal M_k^{\perp}=\{y\in\mathcal H\;\vline\;\forall z\in\mathcal M_k\quad
	\mathfrak I[T_{\rho}(a_kd_k\cdot\lambda)z,y]=0\},\qquad k=1,2,\ldots,n.
\]

Предположим, что справедливо неравенство
\begin{equation}\label{eq3:3.xx}
	\dim\mathcal M-\sum\limits_{k=1}^n\dim\mathcal M_k>n-1.
\end{equation}
Поскольку подпространство
\[
	\tilde{\mathcal M}=\bigoplus_{k=1}^n G_k(\mathcal M_k^{\perp})
\]
имеет в \(\mathcal H\) коразмерность \(\sum_{k=1}^n\dim\mathcal M_k+n-1\), то
подпространство \(\mathcal M\cap\tilde{\mathcal M}\) нетривиально. Следовательно,
найдётся такой набор функций \(y_k\in\mathcal M_k^{\perp}\), где \(k=1,2,\ldots,n\),
что для отвечающей ему функции \(y=\sum_{k=1}^nG_k(y_k)\) будут справедливы соотношения
\[
	\sum\limits_{k=1}^n a_k^{-1}\,\mathfrak I[T_{\rho}(a_kd_k\cdot\lambda)y_k,y_k]=
	\mathfrak I[T_{\rho}(\lambda)y,y]<0.
\]
Тем самым при некотором \(m\in\{1,2,\ldots,n\}\) будет справедливо неравенство \(\mathfrak I
[T_{\rho}(a_md_m\cdot\lambda)y_m,y_m]<0\). Но в этом случае, как легко проверить
непосредственно, подпространство
\[
	\mathcal M'_{m}=\mathcal M_{m}\oplus\operatorname{Lin}\{y_m\}
\]
будет удовлетворять условию
\[
	\forall y\in\mathcal M'_{m}\qquad\mathfrak I[T_{\rho}(a_md_m\cdot\lambda)y,y]
	\leqslant-\varepsilon'_m\|y\|^2_{\mathcal H}
\]
при некотором \(\varepsilon'_m>0\).

Сказанное означает, что, отправляясь от набора подпространств \(\mathcal M_k\), где
\(k=1,2,\ldots,n\), удовлетворяющего условиям~\eqref{eq3:3.xxl} и~\eqref{eq3:3.xx},
всегда можно построить набор подпространств \(\mathcal M'_k\), где \(k=1,2,\ldots,n\),
удовлетворяющий при некоторых \(\varepsilon'_k>0\) условиям
\[
	\forall y\in\mathcal M'_k\qquad\mathfrak I[T_{\rho}(a_k d_k\cdot\lambda)y,y]
	\leqslant-\varepsilon'_k\|y\|^2_{\mathcal H},\qquad k=1,2,\ldots,n,
\]
а также условию
\[
	\sum\limits_{k=1}^n\dim\mathcal M'_k>\sum\limits_{k=1}^n\dim\mathcal M_k.
\]
Отсюда автоматически вытекает справедливость правого неравенства в
соотношении~\eqref{eq3:3.2}. Лемма полностью доказана.
\end{proof}

Утверждение леммы~\ref{lem3:2.1} представляет собой обобщение формулы~(18) из
работы~\cite{SV}.

\begin{proof}[Доказательство теоремы~\ref{tm3:2}]
\textit{Шаг~1.} Введём в рассмотрение непрерывные функции \(\Lambda_{\pm,\varepsilon}\)
вида
\begin{equation}\label{eq4:lambd}
	\forall t\in\mathbb R\qquad\Lambda_{\pm,\varepsilon}(t)=
	e^{-D\nu t/2}\,\varepsilon^{-1}\int\limits_{t}^{t+\varepsilon}
	\ind T_{\rho}\left(\pm e^{\nu \zeta}\right)\,d\zeta.
\end{equation}
Здесь \(\varepsilon\) "--- произвольная положительная постоянная. Из леммы~\ref{lem3:2.1}
следует, что имеют место неравенства
\begin{multline}\label{eq4:taks}
	\forall t\in\mathbb R\qquad\left|\Lambda_{\pm,\varepsilon}(t)-
	\sum\limits_{d_k>0}(a_k\,|d_k|)^{D/2}\Lambda_{\pm,\varepsilon}(t-l_k)-\right.\\
	-\left.\sum\limits_{d_k<0}(a_k\,|d_k|)^{D/2}\Lambda_{\mp,\varepsilon}(t-l_k)\right|
	\leqslant e^{-D\nu t/2}\cdot(n-1),
\end{multline}
где числа \(l_k\) определены соотношением~\eqref{eq:l_k666}. Кроме того, с очевидностью
найдутся такие \(t_0\in\mathbb R\) и \(\varepsilon_0>0\), что будут справедливы тождества
\[
	\forall\varepsilon\in(0,\varepsilon_0)\;\forall t\leqslant t_0\qquad
	\Lambda_{\pm,\varepsilon}(t)=0.
\]
Применяя теперь теоремы~\ref{tmA:1} и~\ref{tmA:2} к функциям
\(\tilde{\Lambda}_{\pm,\varepsilon}\) вида
\[
	\forall t\in\mathbb R\qquad\tilde{\Lambda}_{\pm,\varepsilon}(t)=
	\Lambda_{\pm,\varepsilon}(t-t_0),
\]
находим, что при \(t\to+\infty\) функции \(\Lambda_{\pm,\varepsilon}\) равномерно
по \(\varepsilon\in(0,\varepsilon_0)\) стремятся
к непрерывным \mbox{\(1\)-пе}\-ри\-о\-ди\-че\-ским функциям \(s_{\pm,\varepsilon}\).
При этом в случае, когда среди коэффициентов \(d_k\) имеются отрицательные, справедливо
тождество
\[
	\forall\varepsilon\in(0,\varepsilon_0)\;\forall t\in\mathbb R\qquad
	s_{+,\varepsilon}(t)=s_{-,\varepsilon}(t).
\]

\textit{Шаг~2.} Зафиксируем произвольное \(\delta>0\). Зафиксируем также такое натуральное
число \(N\), что \(e^{-D\nu N/2}<\delta/4\) и
\begin{equation}\label{eq3:3.taks}
	\forall\varepsilon\in(0,\varepsilon_0)\;\forall t\in[N,N+1)\qquad
	|\Lambda_{\pm,\varepsilon}(t)-s_{\pm,\varepsilon}(t)|<\delta/4.
\end{equation}
Из теоремы~\ref{tm2:1} следует, что для некоторого \(\varepsilon_1\in(0,\varepsilon_0)\)
будут справедливы неравенства
\[
	\forall\varepsilon\in(0,\varepsilon_1)\;\forall t\in[N,N+1)\qquad
	|e^{D\nu t/2}\Lambda_{\pm,\varepsilon}(t)-\ind T_{\rho}(\pm e^{\nu t})|
	\leqslant 1.
\]
Отсюда немедленно вытекают неравенства
\[
	\forall\varepsilon,\varepsilon'\in(0,\varepsilon_1)\;\forall t\in[N,N+1)\qquad
	|\Lambda_{\pm,\varepsilon}(t)-\Lambda_{\pm,\varepsilon'}(t)|<\delta/2,
\]
объединяя которые с неравенствами~\eqref{eq3:3.taks}, убеждаемся в справедливости
неравенств
\[
	\forall\varepsilon,\varepsilon'\in(0,\varepsilon_1)\;\forall t\in[N,N+1)\qquad
	|s_{\pm,\varepsilon}(t)-s_{\pm,\varepsilon'}(t)|<\delta.
\]
Суммируя сказанное, получаем, что при \(\varepsilon\searrow 0\) функции
\(s_{\pm,\varepsilon}\) равномерно на \(\mathbb R\) сходятся к некоторым
\mbox{\(1\)-пе}\-ри\-о\-ди\-че\-ским функциям \(s_{\pm}\), для которых справедливы
асимптотики~\eqref{eq3:as}. При этом в случае, когда среди величин \(d_k\)
имеются отрицательные, справедливо тождество~\eqref{eq3:id+}.
Тем самым справедливость первых двух утверждений теоремы доказана.

\textit{Шаг~3.} Для завершения доказательства теоремы остаётся показать справедливость её
утверждения~\eqref{ok:3}. Мы ограничимся рассмотрением случая с функцией \(s_+\) (случай
с функцией \(s_-\) рассматривается аналогично).

На данном шаге будет предполагаться, что все коэффициенты \(d_k\) неотрицательны.

Пусть при некотором \(y\in\mathcal H\) справедливо неравенство \(\mathfrak I[\rho,
|y|^2]>0\). Тогда, согласно определению пучка \(T_{\rho}\), при \(\lambda\gg 0\) будет
справедливо неравенство \(\ind T_{\rho}(\lambda)>0\). Согласно теореме~\ref{tm2:1},
это означает наличие собственных значений пучка \(T_{\rho}\) на положительной
полупрямой.

Обозначим через \(\lambda_1\) наименьшее положительное собственное значение пучка \(T_{\rho}\).
Тогда при любом достаточно малом \(\varepsilon>0\) на полуинтервале \(t\in[\ln\lambda_1/\nu,
\ln\lambda_1/\nu+1-\varepsilon)\) будет выполняться неравенство
\[
	\Lambda_{+,\varepsilon}(t)-\sum\limits_{d_k>0}(a_k\,|d_k|)^{D/2}
	\Lambda_{+,\varepsilon}(t-l_k)\geqslant e^{-D\nu/2}\,\lambda_1^{-1}.
\]
Согласно теореме~\ref{tmA:1} и лемме~\ref{lem3:2.1} это означает, что функция \(s_+\) должна
подчиняться неравенству
\[
	\forall t\in\mathbb R\qquad s_+(t)\geqslant\dfrac{e^{-D\nu/2}}{J\lambda_1},
\]
где положено
\begin{equation}\label{eq3:3.J}
	J=\sum\limits_{d_k\neq 0}l_k\,(a_k\,|d_k|)^{D/2}.
\end{equation}
Следовательно, в рассматриваемом случае функция \(s_+\) положительна.

\textit{Шаг~4.} Пусть теперь среди величин \(d_k\) имеются
отрицательные. Обозначим через \(\lambda_1\) наименьшую из абсолютных величин собственных
значений пучка \(T_{\rho}\). Тогда при любом достаточно малом \(\varepsilon>0\) на полуинтервале
\(t\in[\ln\lambda_1/\nu,\ln\lambda_1/\nu+1-\varepsilon)\) будет выполняться по меньшей
мере одно из двух неравенств
\[
	\Lambda_{\pm,\varepsilon}(t)-\sum\limits_{d_k>0}(a_k\,|d_k|)^{D/2}
	\Lambda_{\pm,\varepsilon}(t-l_k)-\sum\limits_{d_k<0}(a_k\,|d_k|)^{D/2}
	\Lambda_{\mp,\varepsilon}(t-l_k)\geqslant e^{-D\nu/2}\,\lambda_1^{-1}.
\]
Согласно теореме~\ref{tmA:2} и лемме~\ref{lem3:2.1} это означает, что функция \(s_+\) должна
подчиняться неравенству
\[
	\forall t\in\mathbb R\qquad s_+(t)\geqslant\dfrac{e^{-D\nu/2}}{2J\lambda_1},
\]
где величина \(J\) определена соотношением~\eqref{eq3:3.J}. Следовательно, в рассматриваемом
случае функция \(s_+\) также положительна. Тем самым третье утверждение теоремы справедливо.

Теорема полностью доказана.
\end{proof}

\section{Примеры}\label{pt:4}
\subsection{Некоторые конкретные весовые функции} В дальнейшем через \(P_{a,\delta}\),
где \(a\in(0,1/2)\) и \(\delta\in[0,1/3)\), будет обозначаться квадратично
суммируемая самоподобная функция вида
\begin{multline*}
	n=3,\quad a_1=a_2=a,\quad a_2=1-2a,\\ d_1=d_3=1/2+\delta,\quad d_2=-2\delta,
	\quad\beta_1=0,\quad\beta_2=1/2+\delta,\quad\beta_3=1/2-\delta.
\end{multline*}
В частности, \(P_{1/3,0}\) представляет собой хорошо известную канторову лестницу.

Кроме того, через \(\tilde P_{a}\), где \(a\in(0,1/3)\), будет обозначаться функция
\(P_{a,\delta}\), для которой значение параметра \(\delta\) определено условием
\((2-5a)\,\delta=a/2\).

\begin{prop}\label{prop4:1}
Пусть \(D_a\) есть спектральный порядок функции \(\tilde P_a\). Тогда при
\(a\nearrow 1/3\) справедлива асимптотика \(D_a\nearrow 2\).
\end{prop}

\begin{proof}
В рассматриваемом случае имеет место соотношение
\[
	D_a=\dfrac{\ln 9}{\ln (2-5a)-\ln a-\ln(1-2a)}.
\]
Подставляя сюда \(a\nearrow 1/3\), получаем искомое утверждение.
\end{proof}

\subsection{Численные результаты}
\begin{table}[t]
\begin{tabular}{|r|rrr|rrr|c|r|rrr|rrr|}
\cline{1-7} \cline{9-15}
{\(n\)}&\multicolumn{3}{|c|}{\(\lambda_n\)}&\multicolumn{3}{|c|}{\(n/\lambda_n^{\log_6 2}\)}&
\hphantom{!}\hspace{2em}\hphantom{!}&{\(n\)}&\multicolumn{3}{|c|}{\(\lambda_n\)}&
\multicolumn{3}{|c|}{\(n/\lambda_n^{\log_6 2}\)}\\ \cline{1-7} \cline{9-15}
1&\(1,44\cdot 10^1\)&\(\pm\)&\(1\%\)&0,356&\(\pm\)&0,001&
&11&\(2,03\cdot 10^3\)&\(\pm\)&\(1\%\)&0,578&\(\pm\)&0,001\\
2&\(3,53\cdot 10^1\)&\(\pm\)&\(1\%\)&0,504&\(\pm\)&0,001&
&12&\(2,03\cdot 10^3\)&\(\pm\)&\(1\%\)&0,630&\(\pm\)&0,001\\
3&\(1,41\cdot 10^2\)&\(\pm\)&\(1\%\)&0,442&\(\pm\)&0,001&
&13&\(2,27\cdot 10^3\)&\(\pm\)&\(1\%\)&0,654&\(\pm\)&0,001\\
4&\(1,51\cdot 10^2\)&\(\pm\)&\(1\%\)&0,574&\(\pm\)&0,001&
&14&\(2,29\cdot 10^3\)&\(\pm\)&\(1\%\)&0,702&\(\pm\)&0,001\\
5&\(3,26\cdot 10^2\)&\(\pm\)&\(1\%\)&0,533&\(\pm\)&0,001&
&15&\(5,26\cdot 10^3\)&\(\pm\)&\(1\%\)&0,545&\(\pm\)&0,001\\
6&\(3,53\cdot 10^2\)&\(\pm\)&\(1\%\)&0,620&\(\pm\)&0,001&
&16&\(5,26\cdot 10^3\)&\(\pm\)&\(1\%\)&0,582&\(\pm\)&0,001\\
7&\(8,76\cdot 10^2\)&\(\pm\)&\(1\%\)&0,509&\(\pm\)&0,001&
&17&\(9,23\cdot 10^3\)&\(\pm\)&\(1\%\)&0,497&\(\pm\)&0,001\\
8&\(8,76\cdot 10^2\)&\(\pm\)&\(1\%\)&0,582&\(\pm\)&0,001&
&18&\(9,27\cdot 10^3\)&\(\pm\)&\(1\%\)&0,525&\(\pm\)&0,001\\
9&\(1,58\cdot 10^3\)&\(\pm\)&\(1\%\)&0,521&\(\pm\)&0,001&
&19&\(9,59\cdot 10^3\)&\(\pm\)&\(1\%\)&0,547&\(\pm\)&0,001\\
10&\(1,62\cdot 10^3\)&\(\pm\)&\(1\%\)&0,573&\(\pm\)&0,001&
&20&\(9,60\cdot 10^3\)&\(\pm\)&\(1\%\)&0,576&\(\pm\)&0,001\\
\cline{1-7} \cline{9-15}
\end{tabular}

\vspace{0.5cm}
\caption{Оценки первых собственных значений для случая \(a=1/3\), \(\delta=0\)}
\label{tab:1}
\end{table}
В таблице~\ref{tab:1} представлены результаты численных расчётов для первых двадцати
собственных значений задачи Штурма--Лиувилля, весовой функцией в которой выступает
производная функции \(P_{1/3,0}\). Данные таблицы позволяют проиллюстрировать следующее
из теоремы~\ref{tm3:2} утверждение \(\lambda_n\asymp n^{\ln 6/\ln 2}\).

\begin{figure}[t]
\unitlength=0.008cm
\begin{picture}(1000,1000)
\put(0,0){\vector(1,0){1000}}
\put(0,0){\vector(0,1){1000}}
\put(-15,750){\llap{\(0,75\)}}
\put(-15,500){\llap{\(0,5\)}}
\put(-15,1000){\llap{\(1\)}}
\put(15,15){\text{\(0\)}}
\put(1015,15){\text{\(1\)}}
\put(515,15){\text{\(0,5\)}}
\thinlines
\dashline[0]{75}(0,750)(1000,750)
\dashline[0]{75}(0,500)(1000,500)
\path(500,-15)(500,15)
\Thicklines
\path(0,410)(50,400)(100,400)(125,450)(150,500)(200,510)(250,500)
\path(250,500)(300,500)(350,530)(400,630)(450,610)(500,590)(550,560)(600,550)
\path(600,550)(650,530)(700,510)(750,490)(800,490)(850,490)(900,470)(950,460)(1000,460)
\path(0,570)(50,610)(100,650)(125,700)(150,700)(200,720)(250,720)(300,820)
\path(300,820)(350,820)(400,790)(450,770)(500,740)(550,710)(600,690)(650,670)
\path(650,670)(700,640)(750,640)(800,640)(850,620)(900,600)(950,580)(1000,560)
\end{picture}
\caption{Оценки функции \(s_+\) для случая \(a=1/3\), \(\delta=0\)}\label{bil:1}
\end{figure}

Как следует из соотношений~\eqref{eq4:lambd} и~\eqref{eq4:taks}, для функции
\(\rho=P'_{1/3,0}\) и любого положительного \(\lambda\), не являющегося собственным
значением пучка \(T_{\rho}\), справедливы соотношения
\begin{equation}\label{eq5:++-}
	\lambda^{-\log_6 2}\cdot\left(\ind T_{\rho}(\lambda)-2\right)\leqslant
	s_+(\log_6 \lambda)\leqslant\lambda^{-\log_6 2}
	\cdot\left(\ind T_{\rho}(\lambda)+2\right).
\end{equation}
Подставляя в эти соотношения данные таблицы~\ref{tab:1}, устанавливаем справедливость
неравенств
\[
	s_+(\log_6 \lambda_{14}+0)\geqslant 0,60,\qquad
	s_+(\log_6 \lambda_{17}+0)\leqslant 0,56.
\]
Тем самым для рассматриваемой весовой функции \(\rho\) функция \(s_+\) не является
постоянной. Грубые верхняя и нижняя оценки для графика этой функции \(s_+\)
представлены на рисунке~\ref{bil:1}.

Рассмотрим также функцию \(\hat P\), для которой \(\beta_2=2/5\), а прочие параметры
самоподобия совпадают с таковыми для функции \(P_{1/3,0}\). Вес \(\hat\rho=\hat P'\)
является индефинитным, но к нему не может быть применено утверждение~\eqref{ok:2}
теоремы~\ref{tm3:2}. Расчёты показывают, что имеют место равенства
\[
	\ind T_{\hat\rho}(10000)=19,\qquad\ind T_{\hat\rho}(-10000)=3.
\]
Однако для функции \(\rho=\hat\rho\) при любом положительном \(\lambda\), не принадлежащем
спектру \(T_{\rho}\), справедливы соотношения~\eqref{eq5:++-}, а при любом отрицательном
\(\lambda\), не принадлежащем спектру \(T_{\rho}\), справедливы соотношения
\[
	|\lambda|^{-\log_6 2}\cdot\left(\ind T_{\rho}(\lambda)-2\right)\leqslant
	s_-(\log_6 |\lambda|)\leqslant|\lambda|^{-\log_6 2}
	\cdot\left(\ind T_{\rho}(\lambda)+2\right).
\]
Поэтому имеют место неравенства
\[
	s_+(\log_6 10000)\geqslant 0,48,\qquad s_-(\log_6 10000)\leqslant 0,15,
\]
показывающие, что требование наличия отрицательного \(d_k\) существенно для справедливости
утверждения~\eqref{ok:2} теоремы~\ref{tm3:2}.

\vspace{0.3cm}
\textbf{Благодарности:} Авторы благодарят А.~И.~Назарова и А.~А.~Шкаликова за
обсуждение работы и ценные замечания.

\end{document}